\documentclass[11pt]{amsart}
\usepackage[T1]{fontenc}
\usepackage[latin1]{inputenc}
\usepackage{a4wide}
\usepackage{setspace}
\usepackage{xypic}
\usepackage{amssymb}
\usepackage{amsthm}
\usepackage[english]{babel}
\usepackage{latexsym}
\date{}
\newtheorem{thm}{Theorem}[section]
\newtheorem{lem}[thm]{Lemma}

\newtheorem{defn}[thm]{Definition}
\newtheorem*{thm*}{Theorem}
\newtheorem{rem}[thm]{Remark}

\newtheorem{prop}[thm]{Proposition}

\newtheorem{cor}[thm]{Corollary}
\newtheorem{notation}[thm]{Notation}
\newenvironment{f-proof}[1][\sc D\'emonstration.]{\begin{trivlist}
\item[\hskip \labelsep {\bfseries #1}]}{\hfill{$\square$}\end{trivlist}}

\newcommand{\Prod}{\displaystyle\prod}
\newcommand{\fonc}[5]{
 \begin{array}{cccc}
 #1: & #2 & \longrightarrow & #3\\
     & #4 & \longmapsto & #5
 \end{array}
}

\begin{document}

\title[Milne's correcting factor II]{Milne's correcting factor and derived de Rham cohomology II}
\author{Baptiste Morin}

\maketitle

\begin{abstract}
Milne's correcting factor, which appears in the Zeta-value at $s=n$ of a smooth projective variety $X$ over a finite field $\mathbb{F}_q$, is the Euler characteristic of the derived de Rham cohomology of $X/\mathbb{Z}$ modulo the Hodge filtration $F^n$. In this note, we extend this result to arbitrary separated schemes of finite type over $\mathbb{F}_q$ of dimension at most $d$, provided resolution of singularities for schemes of dimension at most $d$ holds. More precisely, we show that Geisser's generalization of Milne's factor, whenever it is well defined, is the Euler characteristic of the $eh$-cohomology with compact support of the derived de Rham complex relative to $\mathbb{Z}$ modulo $F^n$.

\end{abstract}

\footnote{

B. Morin\\
CNRS, IMB\\
Universit\'e de Bordeaux\\
351, cours de la Lib\'eration - F 33405 Talence cedex, France\\ 
e-mail: Baptiste.Morin@math.u-bordeaux.fr\\

The author was supported by ANR-12-BS01-0002 and ANR-12-JS01-0007.}

\section{Introduction}

For any separated scheme $X$ of finite type over the finite field $\mathbb{F}_q$, the special values of the zeta function $Z(X,t):=\prod_{x\in X_0}(1-t^{\mathrm{deg}(x)})^{-1}$ are conjecturally given by 
\begin{equation}\label{LG-Conj}
\mathrm{lim}_{t\rightarrow q^{-n}} Z(X,t)\cdot (1-q^{n}t)^{\rho_n}= \pm \chi(H_{W,c}^*(X,\mathbb{Z}(n)),\cup e)\cdot q^{\chi^{eh}_c(X/\mathbb{F}_q,\mathcal{O},n)}.
\end{equation}
Here  $H_{W,c}^*(X,\mathbb{Z}(n))$ denotes Geisser's "arithmetic cohomology with compact support", $\cup e$ is cup-product with the fundamental class $e\in H^1(W_{\mathbb{F}_q},\mathbb{Z})$ and $q^{\chi^{eh}_c(X/\mathbb{F}_q,\mathcal{O},n)}$ is Geisser's generalization of Milne's correcting factor. The factor $q^{\chi^{eh}_c(X/\mathbb{F}_q,\mathcal{O},n)}$ is well defined under the assumption that resolution of singularities for schemes of dimension $\leq \mathrm{dim}(X)$ holds. The same assumption guaranties that, for $X$ smooth projective, $H_{W,c}^*(X,\mathbb{Z}(n))$ coincides with Weil-\'etale motivic cohomology and $q^{\chi^{eh}_c(X/\mathbb{F}_q,\mathcal{O},n)}$ coincides with Milne's correcting factor. 
For arbitrary $X$, the definitions of  $H_{W,c}^*(X,\mathbb{Z}(n))$ and $q^{\chi^{eh}_c(X/\mathbb{F}_q,\mathcal{O},n)}$  involve $eh$-cohomology with compact support. For instance
$$\chi^{eh}_c(X/\mathbb{F}_q,\mathcal{O},n):=\sum_{i\leq n, j\in\mathbb{Z}}(-1)^{i+j}\cdot (n-i)\cdot \mathrm{dim}_{\mathbb{F}_q}H^j_{c}(X_{eh},\Omega^i)$$
where $H^i_{c}(X_{eh},\Omega^i)$ denotes $eh$-cohomology with compact support of the sheaf of differentials $\Omega^i$. 
Let $\mathbf{Sch}^d/\mathbb{F}_q$ be the category of separated schemes of finite type over $\mathbb{F}_q$ of dimension at most $d$. We say that $R(d)$ holds if any $X\in \mathbf{Sch}^d/\mathbb{F}_q$ admits resolution of singularities (see \cite{Geisser-arith-coh} Definition 2.4 for a precise statement). T. Geisser has shown  in \cite{Geisser-arith-coh} that, if $R(d)$ holds and if the groups $H_{W}^i(Y,\mathbb{Z}(n))$ are finitely generated for any smooth projective variety $Y$ of dimension at most $d$,  then 
$\chi^{eh}_c(X/\mathbb{F}_q,\mathcal{O},n)$ is well defined and (\ref{LG-Conj}) holds for any $X\in \mathbf{Sch}^d/\mathbb{F}_q$.

It was pointed out in \cite{Morin16} that, for $X$ smooth projective,  Milne's correcting factor is the (multiplicative) Euler-Poincar\'e characteristic of the derived de Rham cohomology complex $R\Gamma(X_{Zar},L\Omega^*_{X/\mathbb{Z}}/F^n)$ and that (\ref{LG-Conj}) can be restated in terms of a certain fundamental line. The aim of this note is to show that this remark applies for arbitrary separated schemes of finite type over $\mathbb{F}_q$. More precisely, we denote by $\mathrm{Sh}_{eh}(\mathbf{Sch}^d/\mathbb{F}_q)$ the category of sheaves of sets on the category $\mathbf{Sch}^d/\mathbb{F}_q$ endowed with the $eh$-topology. The resulting  $eh$-topos $\mathrm{Sh}_{eh}(\mathbf{Sch}^d/\mathbb{F}_q)$ is endowed with a structure ring $\mathcal{O}^{eh}$, which is defined as the $eh$-sheafification of the presheaf $X\mapsto \mathcal{O}_X(X)$ on $\mathbf{Sch}^d/\mathbb{F}_q$. We denote by $L\Omega^*_{\mathcal{O}^{eh}/\mathbb{Z}}/F^n$ the derived de Rham complex modulo the Hodge filtration $F^n$ associated with the morphism of ringed topoi
$$(\mathrm{Sh}_{eh}(\mathbf{Sch}^d/\mathbb{F}_q),\mathcal{O}^{eh})\longrightarrow (\mathrm{Spec}(\mathbb{Z}),\mathcal{O}_{\mathrm{Spec}(\mathbb{Z})})$$
where $\mathcal{O}_{\mathrm{Spec}(\mathbb{Z})}$ is the usual structure sheaf on $\mathrm{Spec}(\mathbb{Z})$. 
Then we consider its cohomology with compact support $R\Gamma_c(X_{eh},L\Omega^*_{\mathcal{O}^{eh}/\mathbb{Z}}/F^n)$.
Under the assumption of Theorem \ref{thm}(4) below, one may define the fundamental line
$$\Delta(X/\mathbb{Z},n):=\mathrm{det}_{\mathbb{Z}}R\Gamma_{W,c}(X,\mathbb{Z}(n))\otimes_{\mathbb{Z}} \mathrm{det}_{\mathbb{Z}}R\Gamma_c(X_{eh},L\Omega^*_{\mathcal{O}^{eh}/\mathbb{Z}}/F^n)$$
and its trivialization
$$\lambda_X:\mathbb{R}\stackrel{\sim}{\longrightarrow}\Delta(X/\mathbb{Z},n)\otimes_{\mathbb{Z}}\mathbb{R}$$
which is induced by the acyclic complex
\begin{equation*}
\cdots\stackrel{\cup \theta}{\longrightarrow} H^i_{W,c}(X,\mathbb{Z}(n))_{\mathbb{R}} \stackrel{\cup \theta}{\longrightarrow} H^{i+1}_{W,c}(X,\mathbb{Z}(n))_{\mathbb{R}} \stackrel{\cup \theta}{\longrightarrow}\cdots
\end{equation*}
Here the fundamental class $\theta=\mathrm{Id}_{\mathbb{R}}\in H^1(\mathbb{R},\mathbb{R})="H^1(W_{\mathbb{F}_1},\mathbb{R})"$ is in some sense analogous to $e\in H^1(W_{\mathbb{F}_q},\mathbb{Z})$. We denote by $\zeta^*(X,n)$ the leading coefficient in the Taylor development of $\zeta(X,s)=Z(X,q^{-s})$ near $s=n$.

\begin{thm}\label{thm} Let $X$ be a separated scheme of finite type over $\mathbb{F}_q$ and let $n\in\mathbb{Z}$ be an integer. Assume that $X$ has dimension $d$ and that $R(d)$ holds.
\begin{enumerate}
\item If $X$ is smooth projective, the canonical map
$$R\Gamma(X_{Zar},L\Omega^*_{X/\mathbb{Z}}/F^n)\rightarrow R\Gamma_c(X_{eh},L\Omega^*_{\mathcal{O}^{eh}/\mathbb{Z}}/F^n)$$
is a quasi-isomorphism.
\item The complex $R\Gamma_c(X_{eh},L\Omega^*_{\mathcal{O}^{eh}/\mathbb{Z}}/F^n)$ is bounded with finite cohomology groups.
\item We have
\begin{equation*}\label{part1}\Prod_{i\in\mathbb{Z}} \mid H^i_c(X_{eh},L\Omega^*_{\mathcal{O}^{eh}/\mathbb{Z}}/F^n)\mid^{(-1)^i}\,\,\, = \,\,\, q^{\chi^{eh}_c(X/\mathbb{F}_q,\mathcal{O},n)}.
\end{equation*}
\item Assume moreover that for any smooth projective variety  $Y$ of dimension $\leq d$, the usual Weil-\'etale cohomology groups $H^i_W(Y,\mathbb{Z}(n))$ are finitely generated for all $i$. Then one has 
\begin{eqnarray*}
\Delta(X/\mathbb{Z},n)&=&\mathbb{Z}\cdot\lambda_X\left(\zeta^*(X,n)^{-1}\right).
\end{eqnarray*}
\end{enumerate}
\end{thm}
In particular,  Theorem \ref{thm}(1)--(3) holds (unconditionally) for $\mathrm{dim}(X)\leq 2$ and Theorem \ref{thm}(4) holds for $\mathrm{dim}(X)\leq 1$. This note is organized as follows. We fix some notations and definitions in Section 2. In Section 3, we give the proof of Theorem \ref{thm}, which is based on the following computation of the cohomology sheaves of the complex $L\Lambda_{\mathcal{O}^{eh}}^n L_{\mathcal{O}^{eh}/\mathbb{Z}}$: we define an isomorphism (see Proposition \ref{lem-clef-eh})
$$\mathcal{H}^{i-n}(L\Lambda_{\mathcal{O}^{eh}}^n L_{\mathcal{O}^{eh}/\mathbb{Z}})\simeq \Omega^{i\leq n}_{\mathcal{O}^{eh}/\mathbb{F}_q}$$
where $\Omega^{i\leq n}:=\Omega^{i}$ for $i\leq n$ and $\Omega^{i\leq n}:=0$ for $i>n$. This argument also gives a slightly different proof of the main result of \cite{Morin16}, see Remark \ref{differentproof}.

\section{Preliminaries}

\subsection{The derived de Rham complex}\label{sect-hom-alg}

Given a ring $A$ and an $A$-module $M$, we denote by $\Lambda_A(M)$ (resp. $\Gamma_A(M)$) the exterior $A$-algebra of $M$ (resp. the divided power algebra of $M$, see \cite{Berthelot-Ogus78} App. A), and by $\Lambda^i_A(M)$ (resp. $\Gamma^i_A(M)$) its submodule of homogeneous elements of degree $i$. If $(\mathcal{S},A)$ is a ringed topos and $M$ an $A$-module, one defines $\Lambda_A(M)$, $\Gamma_A(M)$, $\Lambda^i_A(M)$ and $\Gamma^i_A(M)$ as above, internally in $\mathcal{S}$. Then $\Lambda_A(M)$ (resp. $\Gamma_A(M)$) coincides with the sheafification of $U\mapsto \Lambda_{A(U)}(M(U))$ (resp. $U\mapsto \Gamma_{A(U)}(M(U))$). We denote by $L\Lambda^i_A$ the left derived functor of the (non-additive) exterior power functor $\Lambda_A^i$ (see \cite{Illusie71} I.4.2). We often omit the subscript $A$ and simply write $\Lambda^iM$, $\Gamma^iM$  and $L\Lambda^iM$. Let $A\rightarrow B$ be a morphism of rings in $\mathcal{S}$. We denote by $\Omega^1_{B/A}$ the $B$-module of Kähler differentials, we set $\Omega^i_{B/A}:=\Lambda^i_B\Omega^1_{B/A}$ and we denote by $\Omega^{<n}_{B/A}$ the complex of $A$-modules $[\Omega^0_{B/A}\rightarrow \Omega^1_{B/A} \rightarrow \cdots  \rightarrow \Omega^{n-1}_{B/A}]$ put  in degrees $[0,n-1]$. Let $P_A(B)$ be the standard simplicial free resolution of the $A$-algebra $B$ (see \cite{Illusie71} I.1.5.5.6), and let $L_{B/A}$ be the cotangent complex (\cite{Illusie71} II.1). By definition $L_{B/A}$ is the complex of $B$-modules associated with the simplicial $B$-module  $\Omega^{1}_{P_{A}(B)/A}\otimes_{P_{A}(B)} B$. Similarly we define $L\Lambda_B^iL_{B/A}$ as the (actual) complex of $B$-modules associated with the simplicial $B$-module  $\Omega^{i}_{P_{A}(B)/A}\otimes_{P_{A}(B)} B$. 
The derived de Rham complex modulo $F^n$ is defined as the total complex (see \cite{Illusie72} VIII.2.1)
$$L\Omega^{*}_{B/A}/F^n:= \mathrm{Tot}(\Omega^{<n}_{P_{A}(B)/A})$$
which we simply see in this paper as a complex of $A$-modules. The Hodge filtration on $L\Omega^{*}_{B/A}/F^n$ satisfies
$\mathrm{gr}^p(L\Omega^{*}_{B/A}/F^n)\simeq L\Lambda_B^pL_{B/A}[-p]$ for $p<n$ and $\mathrm{gr}^p(L\Omega^{*}_{B/A}/F^n)=0$ otherwise. For example, if $(X,\mathcal{O}_X)$ is a scheme, then $P_{\mathbb{Z}}(\mathcal{O}_X)$ denotes the standard simplicial free resolution of $\mathbb{Z}\rightarrow \mathcal{O}_X$ in the small Zariski topos of the scheme $X$, and $L_{X/\mathbb{Z}}:=L_{\mathcal{O}_X/\mathbb{Z}}$ is the cotangent complex associated with the morphism of schemes $X\rightarrow \mathrm{Spec}(\mathbb{Z})$. 

If $f:\mathcal{S}'\rightarrow \mathcal{S}$ is a morphism of topoi, we write $f^{-1}:\mathcal{S}\rightarrow \mathcal{S}'$ for the set-theoretic inverse image functor of $f$. Let $f:(\mathcal{S}',A')\rightarrow (\mathcal{S},A)$ be a morphism of ringed topoi, i.e. a morphism of topoi $f:\mathcal{S}'\rightarrow \mathcal{S}$ together with a morphism of rings $f^{-1}A\rightarrow A'$ in $\mathcal{S}'$. One defines
$$L\Omega^{*}_{f}/F^n=L\Omega^{*}_{(\mathcal{S}',A')/(\mathcal{S},A)}/F^n:=L\Omega^{*}_{A'/f^{-1}A}/F^n$$
which is a complex of $f^{-1}A$-modules in $\mathcal{S}'$. We denote by $f^*:\mathrm{Mod}(A)\rightarrow \mathrm{Mod}(A')$ the inverse image functor for modules, i.e. $f^*M:=f^{-1}M\otimes_{f^{-1}A}A'$, where $\mathrm{Mod}(A)$ (resp. $\mathrm{Mod}(A')$) is the category of $A$-modules in $\mathcal{S}$ (resp. of $A'$-modules in $\mathcal{S}'$).

\begin{lem}\label{remIll}
Let $f:\mathcal{S}'\rightarrow \mathcal{S}$ be a morphism of topoi and let $A\rightarrow B$ be a morphism of rings in $\mathcal{S}$.  Then we have $f^{-1}(P_A(B))\simeq P_{f^{-1}A}(f^{-1}B)$, $f^{-1}\left(L\Omega^{*}_{B/A}/F^n\right)\simeq L\Omega^*_{f^{-1}B/f^{-1}A}/F^n$, an isomorphism of $f^{-1}B$-modules $f^{-1}(\Omega^i_{B/A})\simeq \Omega^i_{f^{-1}B/f^{-1}A}$ and an isomorphism of complexes of $f^{-1}B$-modules $f^{-1}(L\Lambda_B^iL_{B/A})\simeq L\Lambda_{f^{-1}B}^iL_{f^{-1}B/f^{-1}A}$.
\end{lem}

\begin{proof}
The identifications $f^{-1}(P_A(B))\simeq P_{f^{-1}A}(f^{-1}B)$ and $f^{-1}(\Omega^1_{B/A})\simeq \Omega^1_{f^{-1}B/f^{-1}A}$ follow from the definitions (see \cite{Illusie71} II.1.2.1.4 and \cite{Illusie71} II.1.1.4.1). Moreover we have $f^{-1}(\Lambda^i_{R}(M))\simeq \Lambda^i_{f^{-1}R}(f^{-1}M)$ for any ring $R$ in $\mathcal{S}$ and any $R$-module $M$. The result follows easily.
\end{proof}

\subsection{Derived de Rham cohomology with compact support}

The following definition is due to Thomas Geisser \cite{Geisser-arith-coh}.  Let $\mathbf{Sch}^d/\mathbb{F}_q$ be the category of separated schemes of finite type over $\mathbb{F}_q$ of dimension $\leq d$. 

\begin{defn}
The $eh$-topology on $\mathbf{Sch}^d/\mathbb{F}_q$ is the Grothendieck topology generated by the following coverings: 
\begin{itemize}
\item \'etale coverings
\item  abstract blow-ups: If we have a cartesian square
\[ \xymatrix{
Z'\ar[r]^{i'}\ar[d]^{f'}& X'  \ar[d]^{f}\\
Z\ar[r]^{i}&X  
}
\]
where $f$ is proper, $i$ a closed embedding, and $f$ induces an isomorphism
$X'-Z'\stackrel{\sim}{\longrightarrow} X-Z$, then $(X'\stackrel{f}{\rightarrow} X,Z\stackrel{i}{\rightarrow} X)$ is a covering.
\end{itemize}
\end{defn}
We denote by $\mathrm{PSh}(\mathbf{Sch}^d/\mathbb{F}_q)$ the category of presheaves of sets on $\mathbf{Sch}^d/\mathbb{F}_q$ and by $\mathrm{Sh}_{eh}(\mathbf{Sch}^d/\mathbb{F}_q)$ the topos of $eh$-sheaves of sets on $\mathbf{Sch}^d/\mathbb{F}_q$.
Note that the functor
$$y:\mathbf{Sch}^d/\mathbb{F}_q\hookrightarrow\mathrm{PSh}(\mathbf{Sch}^d/\mathbb{F}_q)\rightarrow \mathrm{Sh}_{eh}(\mathbf{Sch}^d/\mathbb{F}_q),$$
given by composing the Yoneda embedding and $eh$-sheafification, is not fully faithful. Hence the $eh$-topology is not subcanonical. For example, if $X^{\mathrm{red}}$ denotes the maximal reduced closed subscheme of $X\in \mathbf{Sch}^d/\mathbb{F}_q$, then the induced map $y X^{\mathrm{red}}\rightarrow yX$ is an isomorphism. 
If $U$ is an object of $\mathbf{Sch}^d/\mathbb{F}_q$ and $\mathcal{F}$ an  $eh$-sheaf on $\mathbf{Sch}^d/\mathbb{F}_q$, we choose a Nagata compactification
$U\hookrightarrow X$ with closed complement $Z\hookrightarrow X$ (so that $X$ is proper over $\mathbb{F}_q$ and $U$ is open and dense in $X$), and we define
$$R\Gamma_c(U_{eh},\mathcal{F}):=\mathrm{Cone}\left(R\Gamma(X_{eh},\mathcal{F})\rightarrow R\Gamma(Z_{eh},\mathcal{F})\right)[-1].$$
Here $R\Gamma(X_{eh},\mathcal{F})$ denotes the cohomology of the slice topos $\mathrm{Sh}_{eh}(\mathbf{Sch}^d/\mathbb{F}_q)/yX$ with coefficients in $\mathcal{F}\times yX\rightarrow yX$. Equivalently, $R\Gamma(X_{eh},-)$ is the total derived functor of the functor $\mathcal{F}\mapsto\mathcal{F}(X)$. It can be shown that $R\Gamma_c(U_{eh},\mathcal{F})$ does not depend on the compactification (see \cite{Geisser-arith-coh} Proposition 3.2). Then $R\Gamma_c(U_{eh},\mathcal{F})$ is contravariant for proper maps and covariant for open immersions. For an open-closed decomposition $(U\stackrel{j}{\rightarrow}X\stackrel{i}{\leftarrow} Z)$, there is an exact triangle
$$R\Gamma_c(U_{eh},\mathcal{F})\rightarrow R\Gamma_c(X_{eh},\mathcal{F})\rightarrow R\Gamma_c(Z_{eh},\mathcal{F})\rightarrow$$

\begin{notation} The structure ring $\mathcal{O}^{eh}$ on $\mathrm{Sh}_{eh}(\mathbf{Sch}^d/\mathbb{F}_q)$ 
is the $eh$-sheaf associated with the presheaf of rings 
$$\fonc{\mathcal{R}}{(\mathbf{Sch}^d/\mathbb{F}_q)^{op}}{\mathrm{Rings}}{X}{\mathcal{O}_X(X)}.$$
Consider the morphism of ringed topoi
$$\psi:(\mathrm{Sh}_{eh}(\mathbf{Sch}^d/\mathbb{F}_q),\mathcal{O}^{eh})\longrightarrow (\mathrm{Spec}(\mathbb{Z}),\mathcal{O}_{\mathrm{Spec}(\mathbb{Z})})$$
induced by the evident morphism of sites. Let $$L\Omega^*_{\mathcal{O}^{eh}/\mathbb{Z}}/F^n:=L\Omega^*_{\psi}/F^n$$ be the corresponding derived de Rham complex modulo the $n^{\textrm{th}}$-step of the Hodge filtration. Derived de Rham cohomology modulo $F^n$ with compact support is given by
$$X\mapsto R\Gamma_c(X_{eh},L\Omega^*_{\mathcal{O}^{eh}/\mathbb{Z}}/F^n)$$ 
for $X\in \mathbf{Sch}^d/\mathbb{F}_q$. It is covariantly functorial for open immersions and contravariantly functorial for proper maps.
\end{notation}
We now explain our notation $L\Omega^*_{\mathcal{O}^{eh}/\mathbb{Z}}/F^n$. There is a unique morphism of rings $\mathbb{Z}^{eh}\rightarrow \mathcal{O}^{eh}$, where $\mathbb{Z}^{eh}$ denotes the constant sheaf of rings associated with $\mathbb{Z}$ on $\mathrm{Sh}_{eh}(\mathbf{Sch}^d/\mathbb{F}_q)$. Let $L\Omega^*_{\mathcal{O}^{eh}/\mathbb{Z}^{eh}}/F^n$ be the corresponding derived de Rham complex modulo $F^n$. Then we have
\begin{equation}\label{identical}
L\Omega^*_{\psi}/F^n\simeq L\Omega^*_{\mathcal{O}^{eh}/\mathbb{Z}^{eh}}/F^n.
\end{equation}
Indeed, consider the  structure sheaf $\mathcal{O}_{\mathrm{Spec}(\mathbb{Z})}$ and the constant sheaf $\mathbb{Z}$ over the small Zariski topos of $\mathrm{Spec}(\mathbb{Z})$. We have $L_{\mathcal{O}_{\mathrm{Spec}(\mathbb{Z})}/\mathbb{Z}}=0$ (see \cite{Illusie71} II.2.3.1 and II.2.3.6), hence $L_{\mathcal{O}^{eh}/\psi^{-1}\mathcal{O}_{\mathrm{Spec}(\mathbb{Z})}}\simeq L_{\mathcal{O}^{eh}/\psi^{-1}\mathbb{Z}}=L_{\mathcal{O}^{eh}/\mathbb{Z}^{eh}}$. We obtain $L\Lambda^*L_{\mathcal{O}^{eh}/\mathbb{Z}^{eh}}\simeq L\Lambda^*L_{\mathcal{O}^{eh}/\psi^{-1}\mathcal{O}_{\mathrm{Spec}(\mathbb{Z})}}$ hence $$L\Omega^*_{\mathcal{O}^{eh}/\mathbb{Z}^{eh}}/F^n\simeq L\Omega^*_{\mathcal{O}^{eh}/\psi^{-1}\mathcal{O}_{\mathrm{Spec}(\mathbb{Z})}}/F^n:=L\Omega^*_{\psi}/F^n$$
by the Hodge filtration. Finally, we note that $L_{\mathcal{O}^{eh}/\mathbb{Z}}/F^n$ and $L\Omega^*_{\mathcal{O}^{eh}/\mathbb{Z}}/F^n$  could be left-unbounded (see however Corollary \ref{cor-concentration}). 

\subsection{The fundamental line}\label{sectfundline}

For an object $C$ in the derived category of abelian groups such that $H^i(C)$ is finitely generated for all $i$ and $H^i(C)=0$ for almost all $i$, we set
$$\mathrm{det}_{\mathbb{Z}}(C):=\bigotimes_{i\in\mathbb{Z}}\mathrm{det}^{(-1)^i}_{\mathbb{Z}}H^i(C).$$ If  $H^i(C)$ is moreover finite for all $i$, then we call the following isomorphism
$$\mathrm{det}_{\mathbb{Z}}(C)\otimes_{\mathbb{Z}}\mathbb{Q}
\stackrel{\sim}{\rightarrow}\bigotimes_{i\in\mathbb{Z}}\mathrm{det}^{(-1)^i}_{\mathbb{Q}}\left(H^i(C)\otimes_{\mathbb{Z}}\mathbb{Q}\right)\stackrel{\sim}{\rightarrow}\bigotimes_{i\in\mathbb{Z}}\mathrm{det}^{(-1)^i}_{\mathbb{Q}}(0)\stackrel{\sim}{\rightarrow}\mathbb{Q}$$
the canonical $\mathbb{Q}$-trivialization of $\mathrm{det}_{\mathbb{Z}}(C)$. In this situation, the canonical $\mathbb{Q}$-trivialization $\mathrm{det}_{\mathbb{Z}}(C)\otimes_{\mathbb{Z}}\mathbb{Q}\simeq\mathbb{Q}$ identifies $\mathrm{det}_{\mathbb{Z}}(C)$ with $$ \mathbb{Z}\cdot \left(\Prod_{i\in\mathbb{Z}}\vert H^i(C)\vert^{(-1)^{i+1}}\right)\subset \mathbb{Q}.$$
For $X\in \mathbf{Sch}^d/\mathbb{F}_q$, one defines \cite{Geisser-arith-coh}
$$R\Gamma_{W,c}(X,\mathbb{Z}(n)):=R\Gamma(W_{\mathbb{F}_q},R\Gamma_{c}(X_{\overline{\mathbb{F}}_q,eh},\rho^{-1}\mathbb{Z}(n)))$$ 
where $\overline{\mathbb{F}}_q$ is an algebraic closure, $W_{\mathbb{F}_q}$ is the Weil group, $\rho$ is the morphism defined in Lemma \ref{lemrho}, and the $\mathbb{Z}(n)$ on the right hand side is the motivic complex on $\mathbf{Sm}^d/\mathbb{F}_q$. Assuming that $R\Gamma_{W,c}(X,\mathbb{Z}(n))$ and $R\Gamma_c(X_{eh},L\Omega^*_{\mathcal{O}^{eh}/\mathbb{Z}}/F^n)$ are both well defined and perfect, the fundamental line is defined as follows:
$$\Delta(X/\mathbb{Z},n):=\mathrm{det}_{\mathbb{Z}}R\Gamma_{W,c}(X,\mathbb{Z}(n))\otimes_{\mathbb{Z}} \mathrm{det}_{\mathbb{Z}}R\Gamma_c(X_{eh},L\Omega^*_{\mathcal{O}^{eh}/\mathbb{Z}}/F^n).$$
Consider the map $\mathfrak{f}:W_{\mathbb{F}_q}\rightarrow W_{\mathbb{F}_1}:=\mathbb{R}$, and define 
$\theta=\mathrm{Id}_{\mathbb{R}}\in H^1(\mathbb{R},\mathbb{R})$. Then $\mathfrak{f}^*\theta\in H^1(W_{\mathbb{F}_q},\mathbb{R})$ maps the Frobenius $F\in W_{\mathbb{F}_q}$ to $\mathrm{log}(q)\in\mathbb{R}$, whereas $e\in H^1(W_{\mathbb{F}_q},\mathbb{R})$ maps the Frobenius $F$ to $1\in\mathbb{R}$. We have
$$R\Gamma_{W,c}(X,\mathbb{Z}(n))_\mathbb{R}\simeq R\Gamma(W_{\mathbb{F}_q},R\Gamma_{c}(X_{\overline{\mathbb{F}}_q,eh},\rho^{-1}\mathbb{Z}(n))_\mathbb{R}).$$  
So cup-product with the class $\mathfrak{f}^*\theta\in H^1(W_{\mathbb{F}_q},\mathbb{R})$ defines a map
$$H^i_{W,c}(X,\mathbb{Z}(n))_{\mathbb{R}}\stackrel{\cup \theta}{\longrightarrow} H^{i+1}_{W,c}(X,\mathbb{Z}(n))_{\mathbb{R}}$$
which differs from
$$H^i_{W,c}(X,\mathbb{Z}(n))_{\mathbb{R}}\stackrel{\cup e}{\longrightarrow}H^{i+1}_{W,c}(X,\mathbb{Z}(n))_{\mathbb{R}}$$
by the factor $\mathrm{log}(q)$. The complex
\begin{equation*}
\cdots\stackrel{\cup \theta}{\longrightarrow} H^i_{W,c}(X,\mathbb{Z}(n))_{\mathbb{R}} \stackrel{\cup \theta}{\longrightarrow} H^{i+1}_{W,c}(X,\mathbb{Z}(n))_{\mathbb{R}} \stackrel{\cup \theta}{\longrightarrow}\cdots
\end{equation*}
is acyclic \cite{Geisser-arith-coh} hence gives a trivialization
$$\lambda_X:\mathbb{R}\stackrel{\sim}{\longrightarrow}\mathrm{det}_{\mathbb{R}}R\Gamma_{W,c}(X,\mathbb{Z}(n))_{\mathbb{R}}\stackrel{\sim}{\longrightarrow} \Delta(X/\mathbb{Z},n)\otimes_{\mathbb{Z}}\mathbb{R}.$$
where the second isomorphism is induced by the canonical $\mathbb{Q}$-trivialization of $\mathrm{det}_{\mathbb{Z}}R\Gamma_c(X_{eh},L\Omega^*_{\mathcal{O}^{eh}/\mathbb{Z}}/F^n)$, whose existence requires that $R\Gamma_c(X_{eh},L\Omega^*_{\mathcal{O}^{eh}/\mathbb{Z}}/F^n)$ is bounded with finite cohomology groups.

\section{Proof of Theorem \ref{thm}}
We denote by $\mathbf{Sm}^d/\mathbb{F}_q$ the full subcategory of $\mathbf{Sch}^d/\mathbb{F}_q$ consisting of smooth $\mathbb{F}_q$-schemes. We endow $\mathbf{Sm}^d/\mathbb{F}_q$ with the Zariski topology and we denote by $\mathrm{Sh}_{Zar}(\mathbf{Sm}^d/\mathbb{F}_q)$ the corresponding topos.

Recall the following description of the topos $\mathrm{Sh}_{Zar}(\mathbf{Sm}^d/\mathbb{F}_q)$ (see \cite{SGA4} IV.4.10.6). A sheaf $\mathcal{F}$ on $\mathrm{Sh}_{Zar}(\mathbf{Sm}^d/\mathbb{F}_q)$ can be seen as a family of sheaves $\mathcal{F}_X$ on the small Zariski topos $X_{Zar}$ for any $X\in \mathbf{Sm}^d/\mathbb{F}_q$ together with transition maps $\alpha_f:f^{-1}\mathcal{F}_X\rightarrow \mathcal{F}_Y$ for any map $f:Y\rightarrow X$ satisfying $\alpha_{f\circ g}=\alpha_g\circ g^{-1}\alpha_f$ and such that $\alpha_f$ is an isomorphism whenever $f$ is an open immersion. A morphism $\mathcal{F}\rightarrow \mathcal{G}$ in $\mathrm{Sh}_{Zar}(\mathbf{Sm}^d/\mathbb{F}_q)$ is a given by a family of morphisms $\mathcal{F}_X\rightarrow \mathcal{G}_X$ compatible with the transition maps. For any $X\in \mathbf{Sm}^d/\mathbb{F}_q$, the functor
$$\fonc{\mathrm{res}_X}{\mathrm{Sh}_{Zar}(\mathbf{Sm}^d/\mathbb{F}_q)}{X_{Zar}}{\mathcal{F}}{\mathcal{F}_X},$$
mapping the big Zariski sheaf $\mathcal{F}$ to its restriction $\mathcal{F}_X$ to the  small Zariski site of $X$, commutes with arbitrary small limits and colimits. It is therefore the inverse image of a morphism of topoi
$$s_X: X_{Zar}\longrightarrow \mathrm{Sh}_{Zar}(\mathbf{Sm}^d/\mathbb{F}_q)/X \longrightarrow \mathrm{Sh}_{Zar}(\mathbf{Sm}^d/\mathbb{F}_q).$$
In fact the morphism $X_{Zar}\longrightarrow \mathrm{Sh}_{Zar}(\mathbf{Sm}^d/\mathbb{F}_q)/X$ is a section of the morphism
\begin{equation}\label{onemap}
\mathrm{Sh}_{Zar}(\mathbf{Sm}^d/\mathbb{F}_q)/X\simeq \mathrm{Sh}_{Zar}((\mathbf{Sm}^d/\mathbb{F}_q)/X) \longrightarrow X_{Zar}
\end{equation}
which is induced by the evident morphism of sites. The same description of abelian sheaves on $\mathrm{Sh}_{Zar}(\mathbf{Sm}^d/\mathbb{F}_q)$ is valid. We denote by $\mathcal{O}$ the canonical structure ring on $\mathrm{Sh}_{Zar}(\mathbf{Sm}^d/\mathbb{F}_q)$, i.e. $\mathcal{O}(X):=\mathcal{O}_X(X)$ for any $X\in \mathbf{Sm}^d/\mathbb{F}_q$. We have $\mathrm{res}_X(\mathcal{O})=\mathcal{O}_X$ where $\mathcal{O}_X$ denotes the usual structure sheaf on the smooth scheme $X$. As above, a complex of $\mathcal{O}$-modules $\mathcal{F}$ can be seen as family of complexes of $\mathcal{O}_X$-modules $\mathcal{F}_X$ in the small Zariski topos $X_{Zar}$ together with transition maps of complexes of $\mathcal{O}_Y$-modules $\alpha_f:f^*\mathcal{F}_X:=f^{-1}\mathcal{F}_X\otimes_{f^{-1}\mathcal{O}_X}\mathcal{O}_Y\rightarrow \mathcal{F}_Y$ for any map $f:Y\rightarrow X$ satisfying $\alpha_{f\circ g}=\alpha_g\circ g^*\alpha_f$, and such that $\alpha_f$ is an isomorphism whenever $f$ is an open immersion.

We denote by $R(d)$ the condition given in (\cite{Geisser-arith-coh} Definition 2.4). The morphism $\rho$ of the next lemma was defined in (\cite{Geisser-arith-coh} Lemma 2.5), see also (\cite{Suslin-Voevodsky00} Proposition 5.11).

\begin{lem}\label{lemrho} Assume that $R(d)$ holds. Then we have a composite morphism of topoi 
$$\rho:\mathrm{Sh}_{eh}(\mathbf{Sch}^d/\mathbb{F}_q)\stackrel{\sim}{\longrightarrow}\mathrm{Sh}_{eh}(\mathbf{Sm}^d/\mathbb{F}_q)\longrightarrow \mathrm{Sh}_{Zar}(\mathbf{Sm}^d/\mathbb{F}_q)$$
where the first morphism is an equivalence. Moreover we have
\begin{equation}\label{structurering}
\rho^{-1} \mathcal{O}\simeq \mathcal{O}^{eh}.
\end{equation}
\end{lem}
\begin{proof} 
We consider the topology on $\mathbf{Sm}^d/\mathbb{F}_q$ induced by the $eh$-topology on $\mathbf{Sch}^d/\mathbb{F}_q$ (see \cite{SGA4} III.3), and we define $\mathrm{Sh}_{eh}(\mathbf{Sm}^d/\mathbb{F}_q)$ as the topos of sheaves on this site. It follows from $R(d)$ and (\cite{Geisser-arith-coh} Lemma 2.2.b) that $\mathbf{Sm}^d/\mathbb{F}_q$ is a topologically generating full subcategory of $\mathbf{Sch}^d/\mathbb{F}_q$ with respect to the $eh$-topology. By  (\cite{SGA4} III.4.1), the first morphism is an equivalence. The inclusion functor $(\mathbf{Sm}^d/\mathbb{F}_q,Zar)\rightarrow (\mathbf{Sch}^d/\mathbb{F}_q,eh)$ is continuous, i.e. if $\mathcal{F}$ is an $eh$-sheaf on $\mathbf{Sch}^d/\mathbb{F}_q$ then its restriction to $\mathbf{Sm}^d/\mathbb{F}_q$ is a Zariski sheaf. In other words, the induced $eh$-topology on $\mathbf{Sm}^d/\mathbb{F}_q$ is stronger than the Zariski topology; hence the second morphism is well defined.

Let $u:\mathbf{Sm}^d/\mathbb{F}_q\rightarrow \mathbf{Sch}^d/\mathbb{F}_q$ be the inclusion functor, and let  
$f:\mathrm{Sh}_{eh}(\mathbf{Sch}^d/\mathbb{F}_q)\stackrel{\sim}{\rightarrow}\mathrm{Sh}_{eh}(\mathbf{Sm}^d/\mathbb{F}_q)$ be the induced equivalence. We have a commutative square (see \cite{SGA4} III.1.3)
\[ \xymatrix{
\mathrm{Sh}_{eh}(\mathbf{Sch}^d/\mathbb{F}_q)& \mathrm{Sh}_{eh}(\mathbf{Sm}^d/\mathbb{F}_q) \ar[l]_{f^{-1}} \\
\mathrm{PSh}(\mathbf{Sch}^d/\mathbb{F}_q)\ar[u]^{a}&\mathrm{PSh}(\mathbf{Sm}^d/\mathbb{F}_q)  \ar[l]_{u_!} \ar[u]^{a^{\mathrm{sm}}}
}
\]
where the vertical arrows are the associated sheaf functors. Let $\mathcal{F}\in \mathrm{PSh}(\mathbf{Sch}^d/\mathbb{F}_q)$ be a presheaf of sets and let $u^*$ be the right adjoint of $u_!$. The adjunction morphism $u_!u^*\mathcal{F}\rightarrow \mathcal{F}$ is "bicouvrant" (see \cite{SGA4} III.4.1.1) hence $a(u_!u^*\mathcal{F})\stackrel{\sim}{\rightarrow} a(\mathcal{F})$ is an isomorphism (see \cite{SGA4} II.5.3). Since the square above is commutative, we obtain $$f^{-1}\circ a^{\mathrm{sm}}\circ u^*(\mathcal{F})\simeq a\circ u_!\circ u^*(\mathcal{F})\simeq a(\mathcal{F}).$$
So we have an isomorphism of left exact functors $f^{-1}\circ a^{\mathrm{sm}}\circ u^*\simeq a$, hence a similar isomorphism of functors between the categories of ring objects. Let $\mathcal{R}$ (resp. $\mathcal{O}$) be the presheaf of rings on $\mathbf{Sch}^d/\mathbb{F}_q$ (resp. on $\mathbf{Sm}^d/\mathbb{F}_q$) mapping $X$ to $\mathcal{O}_X(X)$. By definition we have $\mathcal{O}=u^*\mathcal{R}$, $\mathcal{O}^{eh}=a(\mathcal{R})$ and $g^{-1}\mathcal{O}=a^{\mathrm{sm}}(\mathcal{O})$, where $g:\mathrm{Sh}_{eh}(\mathbf{Sm}^d/\mathbb{F}_q)\rightarrow \mathrm{Sh}_{Zar}(\mathbf{Sm}^d/\mathbb{F}_q)$ is the morphism of topoi defined above. We obtain
$$\rho^{-1}(\mathcal{O})\simeq f^{-1}\circ g^{-1}(\mathcal{O}) \simeq f^{-1}\circ a^{\mathrm{sm}}(\mathcal{O})\simeq f^{-1}\circ a^{\mathrm{sm}}\circ u^*(\mathcal{R})\simeq a(\mathcal{R})=:\mathcal{O}^{eh}.$$

\end{proof}
We may therefore promote $\rho$ into a morphism of ringed topoi
$$\rho:(\mathrm{Sh}_{eh}(\mathbf{Sch}^d/\mathbb{F}_q),\mathcal{O}^{eh})\longrightarrow (\mathrm{Sh}_{Zar}(\mathbf{Sm}^d/\mathbb{F}_q),\mathcal{O}).$$
For any $X\in\mathbf{Sm}^d/\mathbb{F}_q$, we shall also consider the morphism of ringed topoi obtained by localisation over $X$:
$$\rho_{/X}:(\mathrm{Sh}_{eh}(\mathbf{Sch}^d/\mathbb{F}_q),\mathcal{O}^{eh})/yX\longrightarrow (\mathrm{Sh}_{Zar}(\mathbf{Sm}^d/\mathbb{F}_q),\mathcal{O})/X.$$
We denote by $\mathbb{Z}$ the constant sheaf on either $\mathrm{Sh}_{Zar}(\mathbf{Sm}^d/\mathbb{F}_q)$ or $\mathrm{Sh}_{eh}(\mathbf{Sch}^d/\mathbb{F}_q)$, and we apply the constructions of Section \ref{sect-hom-alg} to the unique  morphism of rings $\mathbb{Z}\rightarrow \mathcal{O}^{eh}$ (respectively $\mathbb{Z}\rightarrow \mathcal{O}$) in the topos $\mathrm{Sh}_{eh}(\mathbf{Sch}^d/\mathbb{F}_q)$ (respectively in the topos $\mathrm{Sh}_{Zar}(\mathbf{Sm}^d/\mathbb{F}_q)$); see (\ref{identical}) and its proof.

\begin{lem}\label{lem-devisse} Assume  $R(d)$. We have
\begin{equation}\label{eq1}
L\Omega^*_{\mathcal{O}^{eh}/\mathbb{Z}}/F^n\simeq \rho^{-1}\left(L\Omega^*_{\mathcal{O}/\mathbb{Z}}/F^n\right)
\end{equation}
and the complex of abelian sheaves $L\Omega^*_{\mathcal{O}/\mathbb{Z}}/F^n$ on $\mathrm{Sh}_{Zar}(\mathbf{Sm}^d/\mathbb{F}_q)$ is given by the complexes of abelian sheaves on $X_{Zar}$
\begin{equation}\label{eq2}
\mathrm{res}_X\left(L\Omega^*_{\mathcal{O}/\mathbb{Z}}/F^n\right)=L\Omega^*_{X/\mathbb{Z}}/F^n:=\mathrm{Tot}(\Omega^{<n}_{P_{\mathbb{Z}}(\mathcal{O}_X)/\mathbb{Z}})
\end{equation}
and obvious transition maps.
Similarly, we have an isomorphism of complexes of $\mathcal{O}^{eh}$-modules
\begin{equation}\label{eq3}
L\Lambda^i_{\mathcal{O}^{eh}}L_{\mathcal{O}^{eh}/\mathbb{Z}}\simeq \rho^{-1}\left( L\Lambda^i_{\mathcal{O}}L_{\mathcal{O}/\mathbb{Z}} \right)
\end{equation}
and the complex of $\mathcal{O}$-modules $L\Lambda^i_{\mathcal{O}}L_{\mathcal{O}/\mathbb{Z}}$ is given by the complexes of $\mathcal{O}_X$-modules
\begin{equation}\label{eq4}
\mathrm{res}_X\left( L\Lambda^i_{\mathcal{O}}L_{\mathcal{O}/\mathbb{Z}} \right)=L\Lambda^i_{\mathcal{O}_X}L_{X/\mathbb{Z}}:=\mathrm{Tot}(\Omega^{i}_{P_{\mathbb{Z}}(\mathcal{O}_X)/\mathbb{Z}}\otimes_{P_{\mathbb{Z}}(\mathcal{O}_X)} \mathcal{O}_X)
\end{equation}
and obvious transition maps. Finally, we have an isomorphism of $\mathcal{O}^{eh}$-modules
\begin{equation}\label{eq5}
\Omega^i_{\mathcal{O}^{eh}/\mathbb{Z}}\simeq \rho^{-1}\Omega^i_{\mathcal{O}/\mathbb{Z}}
\end{equation}
and the $\mathcal{O}$-module $\Omega^i_{\mathcal{O}/\mathbb{Z}}$ is given by the $\mathcal{O}_X$-modules
\begin{equation}\label{eq6}
\mathrm{res}_X\left( \Omega^i_{\mathcal{O}/\mathbb{Z}} \right)=\Omega^i_{X/\mathbb{Z}}
\end{equation}
and obvious transition maps.

\end{lem}

\begin{proof}
The complex $L\Omega^*_{X/\mathbb{Z}}/F^n:= \mathrm{Tot}(\Omega^{<n}_{P_{\mathbb{Z}}(\mathcal{O}_X)/\mathbb{Z}})$ is functorial on the nose in $X\in \mathbf{Sm}^d/\mathbb{F}_q$. Indeed, given a map $f:Y\rightarrow X$, there is a canonical morphism of complexes of abelian sheaves 
\begin{equation}\label{trans1}
f^{-1}L\Omega^*_{X/\mathbb{Z}}/F^n\simeq L\Omega^*_{f^{-1}\mathcal{O}_X/\mathbb{Z}}/F^n\rightarrow L\Omega^*_{Y/\mathbb{Z}}/F^n,
\end{equation} 
where the first map is supplied by  Lemma \ref{remIll} and the second map is induced by the structural morphism $f^{-1}\mathcal{O}_X\rightarrow \mathcal{O}_Y$. The map $f^{-1}L\Omega^*_{X/\mathbb{Z}}/F^n\rightarrow L\Omega^*_{Y/\mathbb{Z}}/F^n$
is an isomorphism of complexes of abelian sheaves if $f:Y\rightarrow X$ is an open immersion.
Similarly, the map $f$ induces a morphism of complexes of $\mathcal{O}_Y$-modules 
\begin{equation}\label{trans2}
f^{*}L\Lambda_{\mathcal{O}_X}^iL_{X/\mathbb{Z}}\simeq L\Lambda_{f^{-1}\mathcal{O}_X}^iL_{f^{-1}\mathcal{O}_X/\mathbb{Z}}\otimes_{f^{-1}\mathcal{O}_X}\mathcal{O}_Y\rightarrow L\Lambda_{\mathcal{O}_Y}^iL_{\mathcal{O}_Y/\mathbb{Z}}
\end{equation}
which is an isomorphism of complexes if $f$ is an open immersion. We apply Lemma \ref{remIll} to the morphism of topoi
$$s_X: X_{Zar}\longrightarrow \mathrm{Sh}_{Zar}(\mathbf{Sm}^d/\mathbb{F}_q)$$ 
and we observe that the transition maps $$f^{-1}\mathrm{res}_X(L\Omega^*_{\mathcal{O}/\mathbb{Z}}/F^n)\rightarrow \mathrm{res}_Y(L\Omega^*_{\mathcal{O}/\mathbb{Z}}/F^n)$$
and
$$f^{*}\mathrm{res}_X(L\Lambda_{\mathcal{O}}^iL_{\mathcal{O}/\mathbb{Z}}/F^n)\rightarrow \mathrm{res}_Y(L\Lambda_{\mathcal{O}}^i L_{\mathcal{O}/\mathbb{Z}}/F^n)$$
can be identified with (\ref{trans1}) and (\ref{trans2}) respectively. This yields (\ref{eq2}) and (\ref{eq4}). We obtain (\ref{eq1})  and (\ref{eq3}) by applying Lemma \ref{remIll} to the morphism
$$\rho: \mathrm{Sh}_{eh}(\mathbf{Sch}^d/\mathbb{F}_q)\longrightarrow \mathrm{Sh}_{Zar}(\mathbf{Sm}^d/\mathbb{F}_q)$$
since we have $\mathcal{O}^{eh}\simeq\rho^{-1}\mathcal{O}$ by (\ref{structurering}). The proof of (\ref{eq5}) and (\ref{eq6}) is similar.
\end{proof}

For a complex $C$ of sheaves of modules on some topos, we denote by $\mathcal{H}^i(C)$ its $i$-th cohomology sheaf.

\begin{lem}\label{lem-clef}
Let $X$ be a smooth separated scheme of finite type over $\mathbb{F}_q$. Then there is
a canonical isomorphism of sheaves of $\mathcal{O}_X$-modules
\begin{equation*}\label{simple}
\mathcal{H}^i(L\Lambda^n L_{X/\mathbb{Z}}[-n])\simeq \Omega^{i\leq n}_{X/\mathbb{F}_q}
\end{equation*}
where $\Omega^{i\leq n}_{X/\mathbb{F}_q}:=\Omega^{i}_{X/\mathbb{F}_q}$ for $0\leq i\leq n$ and $\Omega^{i\leq n}_{X/\mathbb{F}_q}=0$ otherwise. Moreover, for $f:Y\rightarrow X$ a morphism in $\mathbf{Sm}^d/\mathbb{F}_q$, the
square of $\mathcal{O}_Y$-modules 
\[ \xymatrix{
 f^* \mathcal{H}^i(L\Lambda^n L_{X/\mathbb{Z}}[-n])\ar[r]^{\hspace{1cm}\sim}
\ar[d]& f^*\Omega^{i\leq n}_{X/\mathbb{F}_q}  \ar[d]\\
 \mathcal{H}^i(L\Lambda^n L_{Y/\mathbb{Z}}[-n])\ar[r]^{\hspace{1cm}\sim}
& \Omega^{i\leq n}_{Y/\mathbb{F}_q}
}
\]
commutes, where the left vertical map is induced by (\ref{trans2}) and the right vertical map is the evident one.
\end{lem}
\begin{proof}
Let $X$ be a scheme in $\mathbf{Sm}^d/\mathbb{F}_q$. We have
an exact triangle in the derived category $\mathcal{D}(\mathcal{O}_X)$ of $\mathcal{O}_X$-modules (see \cite{Morin16} for details):
$$\mathcal{O}_X[1]\rightarrow L_{X/\mathbb{Z}}\rightarrow \Omega^1_{X/\mathbb{F}_q}\stackrel{\omega_X}{\rightarrow}\mathcal{O}_X[2].$$
Let $U\subset X$ be an affine open subscheme. Then $\omega_U\in \mathrm{Ext}_{\mathcal{O}_U}^2(\Omega^1_{U/\mathbb{F}_q},\mathcal{O}_U)=0$ and there is a unique isomorphism $$\alpha_U:L_{U/\mathbb{Z}}\stackrel{\sim}{\longrightarrow}\mathcal{O}_U[1]\oplus \Omega^1_{U/\mathbb{F}_q}$$ in the derived category $\mathcal{D}(\mathcal{O}_U)$ of $\mathcal{O}_U$-modules, such that $\mathcal{H}^{-1}(\alpha_U):\mathcal{H}^{-1}(L_{U/\mathbb{Z}})\simeq \mathcal{O}_U$ and $\mathcal{H}^{0}(\alpha_U):\mathcal{H}^{0}(L_{U/\mathbb{Z}})\simeq \Omega^1_{X/\mathbb{F}_q}$ are the isomorphisms given by the triangle above. Indeed, the canonical map
$$\mathrm{Hom}_{\mathcal{D}(\mathcal{O}_U)}(L_{U/\mathbb{Z}},\mathcal{O}_U[1]\oplus \Omega^1_{U/\mathbb{F}_q})$$
$$\longrightarrow
 \mathrm{Hom}_{\mathcal{O}_U} (\mathcal{H}^{-1}(L_{U/\mathbb{Z}}),\mathcal{O}_U)\oplus 
 \mathrm{Hom}_{\mathcal{O}_U} (\mathcal{H}^{0}(L_{U/\mathbb{Z}}),\Omega^1_{U/\mathbb{F}_q})$$
is an isomorphism, as it follows from the spectral sequence
$$\Prod_{n\in\mathbb{Z}}\mathrm{Ext}^p (\mathcal{H}^{n}(L_{U/\mathbb{Z}}),\mathcal{H}^{n+q}(\mathcal{O}_U[1]\oplus \Omega^1_{U/\mathbb{F}_q})) \Rightarrow H^{p+q}(R\mathrm{Hom}(L_{U/\mathbb{Z}},\mathcal{O}_U[1]\oplus \Omega^1_{U/\mathbb{F}_q}))$$
and from the fact higher Ext's vanish since $U$ is affine and $\Omega^1_{X/\mathbb{F}_q}$ is locally free of finite rank. Then $\alpha_U$ is functorial in the open affine subscheme $U$, in the sense that, if $V\subseteq U$ is the inclusion of an open affine subscheme $V$, then $\alpha_U\mid V=\alpha_V$ by uniqueness of $\alpha_V$. We obtain the following isomorphism in $\mathcal{D}(\mathcal{O}_U)$ (see \cite{Morin16} for details):
\begin{eqnarray}
& & L\Lambda^n L_{U/\mathbb{Z}}\\
&\simeq & L\Lambda^n([\mathcal{O}_U \stackrel{0}{\rightarrow} \Omega^1_{U/\mathbb{F}_q}][1])\\  
\label{compute}&\simeq &\label{split} [\Gamma^n\mathcal{O}_U\otimes \Lambda^0 \Omega^1_{U/\mathbb{F}_q}\stackrel{0}{\rightarrow} \Gamma^{n-1}\mathcal{O}_U\otimes \Lambda^1\Omega^1_{U/\mathbb{F}_q} \stackrel{0}{\rightarrow}\cdots \stackrel{0}{\rightarrow} \Gamma^0\mathcal{O}_U\otimes \Lambda^n\Omega^1_{U/\mathbb{F}_q}][n]
\end{eqnarray}
where the differential maps are all trivial. This yields a canonical isomorphism of $\mathcal{O}_U$-modules 
$$a_U:\mathcal{H}^i(L\Lambda^n L_{X/\mathbb{Z}}[-n])\mid U\simeq \mathcal{H}^i(L\Lambda^n L_{U/\mathbb{Z}}[-n])\simeq \Gamma^{n-i}\mathcal{O}_U\otimes_{\mathcal{O}_U}\Omega^i_{U/\mathbb{F}_q}$$ for any $i\in\mathbb{Z}$, where $\Gamma^{n-i}\mathcal{O}_U:=0$ for $n-i<0$ and $\Omega^i_{U/\mathbb{F}_q}:=0$ for $i<0$. Moreover, the isomorphisms $a_U$ are compatible with the restriction maps given by inclusions of affine open subsets $V\subseteq U$, in the sense that $(a_U)\mid V=a_V$. 
Covering $X$ by open affine subschemes $U$ (recall that $X$ is separated so that the intersection of two affine open subschemes is affine), the identifications $a_U$ therefore give an isomorphism of sheaves of $\mathcal{O}_X$-modules
\begin{equation*}
\mathcal{H}^i(L\Lambda^n L_{X/\mathbb{Z}}[-n])\simeq \Gamma^{n-i}\mathcal{O}_X\otimes_{\mathcal{O}_X}\Omega^{i}_{X/\mathbb{F}_q}.
\end{equation*}
For $n-i\geq 0$,  the $\mathcal{O}_X$-module $\Gamma^{n-i}\mathcal{O}_X$ is free of rank one with generator $\gamma_{n-i}(1)$ where $1\in\mathcal{O}_X$ is the unit section and $\gamma_{n-i}:\mathcal{O}_X\rightarrow\Gamma^{n-i}\mathcal{O}_X$ the canonical map. So we obtain an isomorphism
\begin{equation}\label{simples}
\mathcal{H}^i(L\Lambda^n L_{X/\mathbb{Z}}[-n])\simeq \Gamma^{n-i}\mathcal{O}_X\otimes_{\mathcal{O}_X}\Omega^{i}_{X/\mathbb{F}_q}\simeq \Omega^{i\leq n}_{X/\mathbb{F}_q}.
\end{equation}

We now check that the isomorphism (\ref{simples}) is functorial in $X\in \mathbf{Sm}^d/\mathbb{F}_q$. Let $Y$ and $X$ be schemes in $\mathbf{Sm}^d/\mathbb{F}_q$ and let $f:Y\rightarrow X$ be an arbitrary map. There is a morphism of exact triangles (see \cite{Illusie71} II.2.1.5)
\[ \xymatrix{
Lf^*\mathcal{O}_X[1]\ar[r]\ar[d]& Lf^* L_{X/\mathbb{Z}}\ar[r]
\ar[d]& Lf^* \Omega^1_{X/\mathbb{F}_q} \ar[r]^{Lf^*\omega_X}
\ar[d]& Lf^*\mathcal{O}_X[2]\ar[d]\\
\mathcal{O}_Y[1]\ar[r]&  L_{Y/\mathbb{Z}}\ar[r]
&  \Omega^1_{Y/\mathbb{F}_q} \ar[r]^{\omega_Y}
& \mathcal{O}_Y[2]
}
\]
Suppose first that $X$ and $Y$ are affine. Then $\omega_X=0$ and $\omega_Y=0$, and the square
\[ \xymatrix{
 Lf^* L_{X/\mathbb{Z}}\ar[r]_{Lf^*\alpha_X\hspace{1cm}}^{\sim \hspace{1cm}}
\ar[d]& f^*\mathcal{O}_X[1]\oplus f^*\Omega^1_{X/\mathbb{F}_q} \ar[d]\\
 L_{Y/\mathbb{Z}}\ar[r]_{\alpha_Y\hspace{0.8cm}}^{\sim \hspace{0.8cm}}& \mathcal{O}_Y[1]\oplus \Omega^1_{Y/\mathbb{F}_q} }
\]
commutes in $\mathcal{D}(\mathcal{O}_Y)$, since a morphism $Lf^* L_{X/\mathbb{Z}}\rightarrow \mathcal{O}_Y[1]\oplus \Omega^1_{Y/\mathbb{F}_q}$ is determined by the morphisms it induces on cohomology (as for $\alpha_X$ above). Hence the bottom square in the following diagram
\[ \xymatrix{
 Lf^* L\Lambda^n L_{X/\mathbb{Z}}\ar[r]^{\sim\hspace{4cm}}
\ar[d]& f^*[\Gamma^{n}_{\mathcal{O}_X}\mathcal{O}_X\otimes\Lambda^{0}_{\mathcal{O}_X}\Omega^1_{X/\mathbb{F}_q} \stackrel{0}{\rightarrow}\cdots \stackrel{0}{\rightarrow} \Gamma^{0}_{\mathcal{O}_X}\mathcal{O}_X\otimes\Lambda^{n}_{\mathcal{O}_X}\Omega^1_{X/\mathbb{F}_q}][n] \ar[d]\\
  L\Lambda^n Lf^* L_{X/\mathbb{Z}}\ar[r]^{\sim\hspace{4.5cm}}
\ar[d]& [\Gamma^{n}_{\mathcal{O}_Y}f^*\mathcal{O}_X\otimes \Lambda^{0}_{\mathcal{O}_Y}f^*\Omega^1_{X/\mathbb{F}_q} \stackrel{0}{\rightarrow}\cdots \stackrel{0}{\rightarrow} \Gamma^{0}_{\mathcal{O}_Y}f^*\mathcal{O}_X\otimes \Lambda^{n}_{\mathcal{O}_Y}f^*\Omega^1_{X/\mathbb{F}_q}][n] \ar[d]\\
 L\Lambda^n L_{Y/\mathbb{Z}}\ar[r]^{\sim \hspace{3cm}}
& [\Gamma^{n}\mathcal{O}_Y\otimes\Omega^0_{Y/\mathbb{F}_q} \stackrel{0}{\rightarrow}\cdots \stackrel{0}{\rightarrow}  \Gamma^{0}\mathcal{O}_Y\otimes\Omega^n_{Y/\mathbb{F}_q}][n]
}
\]
commutes as well (see \cite{Illusie71} I.4.3.1.3). Here the top left vertical map is induced by the derived  version $Lf^*L\Lambda^n_{\mathcal{O}_X}\rightarrow L\Lambda^n_{\mathcal{O}_Y}Lf^*$ of the natural transformation  $f^*\Lambda^n_{\mathcal{O}_X}\rightarrow \Lambda^n_{\mathcal{O}_Y}f^*$, and the top right vertical map is induced by $f^*\Lambda^i_{\mathcal{O}_X}\rightarrow \Lambda^i_{\mathcal{O}_Y}f^*$ and $f^*\Gamma^{n-i}_{\mathcal{O}_X}\rightarrow \Gamma^{n-i}_{\mathcal{O}_Y}f^*$. It follows that the upper square in the previous diagram commutes. 
Since the cohomology sheaves of $L\Lambda^n L_{X/\mathbb{Z}}$ are flat $\mathcal{O}_X$-modules, we have  the isomorphism
$$f^* \mathcal{H}^i(L\Lambda^n L_{X/\mathbb{Z}}[-n])\stackrel{\sim}{\rightarrow} \mathcal{H}^i(Lf^*L\Lambda^n L_{X/\mathbb{Z}}[-n]).$$
We obtain the following commutative square of $\mathcal{O}_Y$-modules
\[ \xymatrix{
 f^* \mathcal{H}^i(L\Lambda^n L_{X/\mathbb{Z}}[-n])\ar[r]^{\hspace{-0.2cm}\sim}
\ar[d]& f^*(\Gamma^{n-i}\mathcal{O}_X\otimes_{\mathcal{O}_X}\Omega^{i}_{X/\mathbb{F}_q})  \ar[d]\\
 \mathcal{H}^i(L\Lambda^n L_{Y/\mathbb{Z}}[-n])\ar[r]^{\hspace{0cm}\sim}
& \Gamma^{n-i}\mathcal{O}_Y\otimes_{\mathcal{O}_Y}\Omega^{i}_{Y/\mathbb{F}_q}
}
\]
where $X$ and $Y$ are affine schemes in $\mathbf{Sm}^d/\mathbb{F}_q$. Let $f:Y\rightarrow X$ be a map between arbitrary $X$, $Y$ in $\mathbf{Sm}^d/\mathbb{F}_q$. Covering $Y$ and $X$ by affine  open subschemes (compatibly with $f$) we see that the previous square commutes for arbitrary $X$ and $Y$.  The result follows because the identification of $\mathcal{O}_X$-modules $\Gamma^{n-i}_{\mathcal{O}_X}\mathcal{O}_X\simeq \mathcal{O}_X$ is functorial in $X$. Indeed, the map $f^*(\Gamma_{\mathcal{O}_X}^{n-i}\mathcal{O}_X)\rightarrow \Gamma^{n-i}_{\mathcal{O}_Y}\mathcal{O}_Y$ maps $\gamma_{n-i}(1)$ to itself.
\end{proof}
\begin{rem}\label{onerem}
An isomorphism of the form $$L\Lambda^n L_{X/\mathbb{Z}}\simeq [\mathcal{O}_X\stackrel{0}{\rightarrow} \Omega^1_{X/\mathbb{F}_q} \stackrel{0}{\rightarrow}\cdots \stackrel{0}{\rightarrow} \Omega^n_{X/\mathbb{F}_q}][n]$$ is false in general, e.g. take $n=1$ and $X$ such that $\alpha_X\neq 0$ (i.e. such that $X$ has no lifting over $\mathbb{Z}/p^2\mathbb{Z}$).
\end{rem}
\begin{rem}\label{differentproof}
In order to prove the main result of \cite{Morin16}, one may use Lemma \ref{lem-clef} above instead of (\cite{Morin16} Lemma 2). 
\end{rem}

\begin{prop}\label{lem-clef-eh} Assume  $R(d)$. There is a canonical isomorphism of sheaves of $\mathcal{O}^{eh}$-modules
\begin{equation*}\label{simple}
\mathcal{H}^i(L\Lambda^n L_{\mathcal{O}^{eh}/\mathbb{Z}}[-n])\simeq \Omega^{i\leq n}_{\mathcal{O}^{eh}/\mathbb{F}_q}
\end{equation*}
where $\Omega^{i\leq n}_{\mathcal{O}^{eh}/\mathbb{F}_q}= \Omega^{i}_{\mathcal{O}^{eh}/\mathbb{F}_q}$ for $0\leq i\leq n$ and $\Omega^{i\leq n}_{\mathcal{O}^{eh}/\mathbb{F}_q}=0$ otherwise.
\end{prop}
\begin{proof}
We first work in the ringed topos $(\mathrm{Sh}_{Zar}(\mathbf{Sm}^d/\mathbb{F}_q),\mathcal{O})$. By exactness of $\mathrm{res}_X$, Lemma \ref{lem-devisse}(\ref{eq4})
and Lemma \ref{lem-clef}, we have
\begin{eqnarray*}
\mathrm{res}_X(\mathcal{H}^i(L\Lambda^n L_{\mathcal{O}/\mathbb{Z}}[-n]))&\simeq & \mathcal{H}^i(\mathrm{res}_X(L\Lambda^n L_{\mathcal{O}/\mathbb{Z}}[-n]))\\
&\simeq & \mathcal{H}^i(L\Lambda^n L_{X/\mathbb{Z}}[-n]))\\
&\simeq & \Omega^{i\leq n}_{X/\mathbb{F}_q}
\end{eqnarray*}
for any $X$ in $\mathbf{Sm}^d/\mathbb{F}_q$. Moreover, for a morphism $f:Y\rightarrow X$ in $\mathbf{Sm}^d/\mathbb{F}_q$, the transition map
$$\alpha_f:f^*\mathrm{res}_X(\mathcal{H}^i(L\Lambda^n L_{\mathcal{O}/\mathbb{Z}}[-n]))\longrightarrow \mathrm{res}_Y(\mathcal{H}^i(L\Lambda^n L_{\mathcal{O}/\mathbb{Z}}[-n]))$$
may be identified with
the canonical map (see Lemma \ref{lem-devisse})
$$f^*\mathcal{H}^i(L\Lambda^n L_{X/\mathbb{Z}}[-n])\longrightarrow \mathcal{H}^i(L\Lambda^n L_{Y/\mathbb{Z}}[-n]) $$
which in turn may be identified with
the canonical map
$$f^*\Omega^{i\leq n}_{X/\mathbb{F}_q}\longrightarrow \Omega^{i\leq n}_{Y/\mathbb{F}_q}$$
by Lemma \ref{lem-clef}.
In view of (\ref{eq6}), we obtain an isomorphism
\begin{equation}\label{equa}
\mathcal{H}^i(L\Lambda^n L_{\mathcal{O}/\mathbb{Z}}[-n])\simeq \Omega^{i\leq n}_{\mathcal{O}/\mathbb{F}_q}
\end{equation}
of $\mathcal{O}$-modules in the topos $\mathrm{Sh}_{Zar}(\mathbf{Sm}^d/\mathbb{F}_q)$.
By Lemma \ref{lem-devisse} (\ref{eq3}) and by exactness of $\rho^{-1}$, we have 
\begin{equation}\label{eqa}
\mathcal{H}^i(L\Lambda^n L_{\mathcal{O}^{eh}/\mathbb{Z}})\simeq \mathcal{H}^i(\rho^{-1} L\Lambda^n L_{\mathcal{O}/\mathbb{Z}})\simeq \rho^{-1} \mathcal{H}^i(L\Lambda^n L_{\mathcal{O}/\mathbb{Z}}).
\end{equation}
By (\ref{eqa}), (\ref{equa}) and (\ref{eq5}), we obtain
$$\mathcal{H}^i(L\Lambda^n L_{\mathcal{O}^{eh}/\mathbb{Z}}[-n])\simeq \rho^{-1} \mathcal{H}^i(L\Lambda^n L_{\mathcal{O}/\mathbb{Z}}[-n])\simeq \rho^{-1} \Omega^{i\leq n}_{\mathcal{O}/\mathbb{F}_q}\simeq \Omega^{i\leq n}_{\mathcal{O}^{eh}/\mathbb{F}_q}.$$
\end{proof}
\begin{rem}
One may think of trying to prove Proposition \ref{lem-clef-eh} more directly using the exact triangle
$$\mathcal{O}^{eh}[1]\rightarrow L_{\mathcal{O}^{eh}/\mathbb{Z}}\rightarrow \Omega^1_{\mathcal{O}^{eh}/\mathbb{F}_q}\stackrel{\omega^{eh}}{\rightarrow}\mathcal{O}^{eh}[2],$$
which is the image by $\rho^{-1}$ of the exact triangle
$$\mathcal{O}[1]\rightarrow L_{\mathcal{O}/\mathbb{Z}}\rightarrow \Omega^1_{\mathcal{O}/\mathbb{F}_q}\stackrel{\omega}{\rightarrow}\mathcal{O}[2]$$
in the derived category of $\mathcal{O}$-modules on $\mathrm{Sh}_{Zar}(\mathbf{Sm}^d/\mathbb{F}_q)$. A direct computation of $L\Lambda^n L_{\mathcal{O}^{eh}/\mathbb{Z}}$ as in (\ref{compute}) would not work since the extension $\omega$ is non-trivial by Remark \ref{onerem} and Lemma \ref{lem-devisse}.
\end{rem}
The following corollary follows immediately from Proposition \ref{lem-clef-eh}.
\begin{cor}\label{cor-concentration}
If $R(d)$ holds then $L\Lambda^p_{\mathcal{O}^{eh}}L_{\mathcal{O}^{eh}/\mathbb{Z}}$ is concentrated in degrees $[-p,0]$ and $L\Omega^*_{\mathcal{O}^{eh}/\mathbb{Z}}/F^n$ is concentrated in degrees $[0,n-1]$.
\end{cor}

\begin{cor}\label{cor-comparison}
Let $X$ be a smooth projective scheme over $\mathbb{F}_q$ of dimension $d$ and let $n\in\mathbb{Z}$ be an integer. If  $R(d)$ holds then the canonical maps
$$R\Gamma(X_{Zar},L\Lambda^{p}_{\mathcal{O_X}}L_{X/\mathbb{Z}})\rightarrow R\Gamma(X_{eh},L\Lambda^{p}_{\mathcal{O}^{eh}}L_{\mathcal{O}^{eh}/\mathbb{Z}})$$
and
$$R\Gamma(X_{Zar},L\Omega^*_{X/\mathbb{Z}}/F^n)\rightarrow R\Gamma(X_{eh},L\Omega^*_{\mathcal{O}^{eh}/\mathbb{Z}}/F^n)$$
are quasi-isomorphisms.
\end{cor}
\begin{proof}
The morphism of ringed topoi
\begin{equation}\label{morphismofringedtopoi}
(\mathrm{Sh}_{eh}(\mathbf{Sch}^d/\mathbb{F}_q),\mathcal{O}^{eh})/y(X)\stackrel{\rho_{/X}}{\longrightarrow} (\mathrm{Sh}_{Zar}(\mathbf{Sm}^d/\mathbb{F}_q),\mathcal{O})/X \stackrel{(\ref{onemap})}{\longrightarrow} (X,\mathcal{O}_X).
\end{equation}
induces a morphism of (derived Hodge to de Rham) spectral sequences 
from
\begin{equation}\label{HtoDR-Zar}
E_1^{p,q}=H^q(X_{Zar},L\Lambda^{p<n}L_{X/\mathbb{Z}})\Longrightarrow H^{p+q}(X_{Zar},L\Omega^*_{X/\mathbb{Z}}/F^n)
\end{equation}
to
\begin{equation}\label{HtoDR-eh}
'E_1^{p,q}=H^q(X_{eh},L\Lambda^{p<n}L_{\mathcal{O}^{eh}/\mathbb{Z}})\Longrightarrow H^{p+q}(X_{eh},L\Omega^*_{\mathcal{O}^{eh}/\mathbb{Z}}/F^n).
\end{equation}
Here (\ref{HtoDR-Zar}) and (\ref{HtoDR-eh}) are obtained  (using Corollary \ref{cor-concentration}) as spectral sequences for the hypercohomology of filtered complexes. 
One is therefore reduced to showing that the maps
$$H^q(X_{Zar},L\Lambda^{p}L_{X/\mathbb{Z}})\rightarrow H^q(X_{eh},L\Lambda^{p}L_{\mathcal{O}^{eh}/\mathbb{Z}})$$
are isomorphisms. By Lemma \ref{lem-clef}, Proposition \ref{lem-clef-eh} and Corollary \ref{cor-concentration}, the morphism (\ref{morphismofringedtopoi}) induces a morphism of hypercohomology spectral sequences
from
\begin{equation*}\label{HyperSS-Zar}
E_2^{i,j}=H^i(X_{Zar},\Omega^{j\leq p}_{X/\mathbb{F}_q})\Longrightarrow H^{i+j}(X_{Zar},L\Lambda^{p}L_{X/\mathbb{Z}}[-p])
\end{equation*}
to
\begin{equation*}\label{HyperSS-eh}
'E_2^{i,j}=H^i(X_{eh},\Omega^{j\leq p}_{\mathcal{O}^{eh}/\mathbb{F}_q})\Longrightarrow H^{i+j}(X_{eh},L\Lambda^{p}L_{\mathcal{O}^{eh}/\mathbb{Z}}[-p]).
\end{equation*}
One is therefore reduced to showing that the map 
$$H^i(X_{Zar},\Omega^j_{X/\mathbb{F}_q})\rightarrow H^i(X_{eh},\Omega^j_{\mathcal{O}^{eh}/\mathbb{F}_q})$$
is an isomorphism for any $i,j$. Assuming $R(d)$, this follows from (\cite{Geisser-arith-coh} Theorem 4.7) since $\Omega^j_{\mathcal{O}^{eh}/\mathbb{F}_q}\simeq \rho^{-1}\Omega^j_{\mathcal{O}/\mathbb{F}_q}$.
\end{proof}
Recall from the introduction that one defines
$$\chi^{eh}_c(X/\mathbb{F}_q,\mathcal{O},n):=\sum_{i\leq n, j\in\mathbb{Z}}(-1)^{i+j}\cdot (n-i)\cdot \mathrm{dim}_{\mathbb{F}_q}H^j_{c}(X_{eh},\Omega_{\mathcal{O}^{eh}/\mathbb{F}_q}^i).$$

\begin{cor} Let $X$ be a separated scheme of finite type over $\mathbb{F}_q$ of dimension $d$ and let $n\in\mathbb{Z}$ be an integer. If $R(d)$ holds then the complex $R\Gamma_c(X_{eh},L\Omega^*_{\mathcal{O}^{eh}/\mathbb{Z}}/F^n)$ is bounded with finite cohomology groups, and we have 
\begin{equation*}\label{part1}\Prod_{i\in\mathbb{Z}} \mid H^i_c(X_{eh},L\Omega^*_{\mathcal{O}^{eh}/\mathbb{Z}}/F^n)\mid^{(-1)^i}\,\,\, = \,\,\, q^{\chi^{eh}_c(X/\mathbb{F}_q,\mathcal{O},n)}.
\end{equation*}
\end{cor}
\begin{proof}
We consider the spectral sequences
\begin{equation}\label{HtoDR-eh-c}
H_c^q(X_{eh},L\Lambda^{p<n}L_{\mathcal{O}^{eh}/\mathbb{Z}})\Longrightarrow H_c^{p+q}(X_{eh},L\Omega^*_{\mathcal{O}^{eh}/\mathbb{Z}}/F^n)
\end{equation}
and
\begin{equation}\label{HyperSS-eh-c}
H_c^i(X_{eh},\Omega^{j\leq p}_{\mathcal{O}^{eh}/\mathbb{F}_q})\Longrightarrow H_c^{i+j}(X_{eh},L\Lambda^{p}L_{\mathcal{O}^{eh}/\mathbb{Z}}[-p]).
\end{equation}
In view of Corollary \ref{cor-concentration} and the isomorphism (see \cite{Geisser-arith-coh} Remark before Lemma 3.5)
$$R\Gamma_c(X_{eh},-)\simeq R\mathrm{Hom}(\mathbb{Z}_{eh}^{c}(X),-)$$
(\ref{HtoDR-eh-c}) and (\ref{HyperSS-eh-c}) may be obtained as spectral sequences for the hypercohomology  of filtered complexes. The complex $R\Gamma_c(X_{eh},\Omega^j_{\mathcal{O}^{eh}/\mathbb{F}_q})\simeq R\Gamma_c(X_{eh},\rho^{-1}\Omega^j_{\mathcal{O}/\mathbb{F}_q})$ is bounded with finite cohomology groups by (\cite{Geisser-arith-coh} Corollary 4.8). In view of (\ref{HtoDR-eh-c}) and (\ref{HyperSS-eh-c}), the complexes $R\Gamma_c(X_{eh},L\Lambda^pL_{\mathcal{O}^{eh}/\mathbb{Z}}/F^n)$ and $R\Gamma_c(X_{eh},L\Omega^*_{\mathcal{O}^{eh}/\mathbb{Z}}/F^n)$ are also  bounded with finite cohomology groups. By (\cite{Morin16} Lemma 1), the spectral sequences (\ref{HtoDR-eh-c}) and (\ref{HyperSS-eh-c}) give isomorphisms
\begin{eqnarray}
&&\mathrm{det}_{\mathbb{Z}}R\Gamma_c(X_{eh},L\Omega^*_{\mathcal{O}^{eh}/\mathbb{Z}}/F^n)\\
&\stackrel{\sim}{\longrightarrow}&\bigotimes_{p<n}\mathrm{det}^{(-1)^{p}}_{\mathbb{Z}}R\Gamma_c(X_{eh},L\Lambda^pL_{\mathcal{O}^{eh}/\mathbb{Z}})\\
&\stackrel{\sim}{\longrightarrow}&\bigotimes_{p<n}\mathrm{det}_{\mathbb{Z}}R\Gamma_c(X_{eh},L\Lambda^pL_{\mathcal{O}^{eh}/\mathbb{Z}}[-p])\\
\label{last}&\stackrel{\sim}{\longrightarrow}&\bigotimes_{p<n}\left(\bigotimes_{i\leq p,j}\mathrm{det}_{\mathbb{Z}}^{(-1)^{i+j}} H^{j}_c(X_{eh},\Omega^i_{\mathcal{O}^{eh}/\mathbb{F}_p})\right)
\end{eqnarray}
such that the square of isomorphisms
\[ \xymatrix{
\left(\mathrm{det}_{\mathbb{Z}}R\Gamma_c(X_{eh},L\Omega^*_{\mathcal{O}^{eh}/\mathbb{Z}}/F^n)\right)_{\mathbb{Q}}\ar[r]\ar[d]^{\gamma}&
\ar[d]^{\gamma'}\left(\bigotimes_{p<n}\bigotimes_{i\leq p,j}\mathrm{det}^{(-1)^{i+j}}_{\mathbb{Z}}H_c^{j}(X_{eh},\Omega^i_{\mathcal{O}^{eh}/\mathbb{F}_q})\right)_{\mathbb{Q}}\\
\mathbb{Q}\ar[r]^{\mathrm{Id}}&\mathbb{Q}
}
\]
commutes, where the top horizontal map is induced by (\ref{last}), and the vertical isomorphisms are the canonical trivializations (see Section \ref{sectfundline}). The result follows:
\begin{eqnarray*}
&&\mathbb{Z}\cdot \left(\Prod_{i\in\mathbb{Z}} \mid H_c^i(X_{eh},L\Omega^*_{\mathcal{O}^{eh}/\mathbb{Z}}/F^n)\mid^{(-1)^i}\right)^{-1}\\
&=&\gamma\left(\mathrm{det}_{\mathbb{Z}}R\Gamma_c(X_{eh},L\Omega^*_{\mathcal{O}^{eh}/\mathbb{Z}}/F^n)\right)\\
&=&\gamma'\left(\bigotimes_{p<n} \bigotimes_{i\leq p,j}\mathrm{det}^{(-1)^{i+j}}_{\mathbb{Z}}H_c^{j}(X_{eh},\Omega^i_{\mathcal{O}^{eh}/\mathbb{F}_q})\right)\\
&=&\mathbb{Z}\cdot q^{-\chi^{eh}_c(X/\mathbb{F}_q,\mathcal{O},n)}.
\end{eqnarray*}\\
\end{proof}
Recall from Section \ref{sectfundline} the definitions of $\Delta(X/\mathbb{Z},n)$ and $\lambda_X$.

\begin{cor}
Let $X$ be a separated scheme of finite type over $\mathbb{F}_q$ of dimension $d$ and let $n\in\mathbb{Z}$ be an integer. Assume that for any smooth projective variety  $Y$ of dimension $\leq d$ the Weil-\'etale cohomology groups $H^i_W(Y,\mathbb{Z}(n))$ are finitely generated for all $i$. If $R(d)$ holds, then one has 
\begin{eqnarray*}
\Delta(X/\mathbb{Z},n)&=&\mathbb{Z}\cdot\lambda_X\left(\zeta^*(X,n)^{-1}\right).
\end{eqnarray*}
\end{cor}
\begin{proof}
All the schemes we consider in this proof are in $\mathbf{Sch}^d/\mathbb{F}_q$. For an open-closed decomposition $(U\stackrel{j}{\rightarrow}X\stackrel{i}{\leftarrow} Z)$, we have exact triangles
$$R\Gamma_c(U_{eh},L\Omega^*_{\mathcal{O}^{eh}/\mathbb{Z}}/F^n)\rightarrow R\Gamma_c(X_{eh},L\Omega^*_{\mathcal{O}^{eh}/\mathbb{Z}}/F^n)\rightarrow R\Gamma_c(Z_{eh},L\Omega^*_{\mathcal{O}^{eh}/\mathbb{Z}}/F^n)$$
and
\begin{equation}\label{Wet-triangle}
R\Gamma_{W,c}(U,\mathbb{Z}(n))\rightarrow R\Gamma_{W,c}(X,\mathbb{Z}(n))\rightarrow R\Gamma_{W,c}(Z,\mathbb{Z}(n)).
\end{equation}
Moreover, the triangle (\ref{Wet-triangle}) is compatible (in the obvious sense) with $\cup\theta$. This gives an isomorphism
\begin{equation}\label{det-triangle}
\Delta(X/\mathbb{Z},n)\simeq \Delta(U/\mathbb{Z},n)\otimes_{\mathbb{Z}} \Delta(Z/\mathbb{Z},n)
\end{equation}
such that the square of isomorphisms
\[ \xymatrix{
\mathbb{R}\ar[r]^{\mathrm{Id}}\ar[d]^{\lambda_X}& \mathbb{R}
\ar[d]^{\lambda_U\otimes \lambda_Z}\\
\Delta(X/\mathbb{Z},n)_{\mathbb{R}}\ar[r]^{(\ref{det-triangle})_\mathbb{R}\hspace{1.5cm}}&\Delta(U/\mathbb{Z},n)_{\mathbb{R}}\otimes_{\mathbb{R}} \Delta(Z/\mathbb{Z},n)_{\mathbb{R}}
}
\]
commutes. Similarly, one has
$$\zeta^*(X,n)=\zeta^*(U,n)\cdot \zeta^*(Z,n).$$
It follows that if the result is true for two out of the three schemes $(X,U,Z)$ then it is true for the third. Moreover, the result is true for $X$ smooth projective by \cite{Morin16}, Corollary \ref{cor-comparison} and (\cite{Geisser-arith-coh} Theorem 4.3). It follows for arbitrary $X\in\mathbf{Sch}^d/\mathbb{F}_q$ by (\cite{Geisser-arith-coh} Lemma 2.7).
\end{proof}

\end{document}